\documentclass[10pt]{amsart}
\usepackage{amsfonts,amssymb}
\def\zz{\mathbb{Z}}

\def\cc{\mathbb{C}}
\def\pp{\mathbb{P}}

\def\hh{\mathbb{H}}

\def\isom{\, {\buildrel \sim \over \longrightarrow} \,}

\def\div{\mathrm{div} \, }

\theoremstyle{plain}
\newtheorem{thm}{Theorem}[section]
\newtheorem{lem}[thm]{Lemma}
\newtheorem{prop}[thm]{Proposition}
\newtheorem{cor}[thm]{Corollary}

\theoremstyle{definition}
\newtheorem{rem}[thm]{Remark}
\newtheorem{example}[thm]{Example}
\newtheorem{definition}[thm]{Definition}

\numberwithin{equation}{section}

\begin{document}

\title{		Theta functions on the theta divisor}
\author{	Robin de Jong}
%\address{	University of Leiden \\
%		The Netherlands}
%\email{		rdejong@math.leidenuniv.nl}
%\subjclass{Primary 14K25; secondary: 14H42, 14H55, 14K20}

\thanks{The author acknowledges a grant from the Netherlands Organisation for
Scientific Research (NWO)}

\begin{abstract}
We show that the gradient and the hessian  
of the Riemann theta function in dimension $n$
can be combined to give a theta function of order
$n+1$ and modular weight $(n+5)/2$ defined on the theta divisor. 
It can be seen that the zero locus of this theta function essentially gives
the ramification locus of the Gauss map. 
For jacobians this leads to a description in terms of theta functions 
and their derivatives of the Weierstrass point locus on the associated Riemann surface.
\end{abstract}

\maketitle

\thispagestyle{empty} 

% introduction

In the analytic theory of the Riemann theta function a natural place is taken by the
study of the first and second order terms of its Taylor series expansion along
the theta divisor. The first order term essentially gives the gradient, and
hence the tangent bundle on the smooth locus, whereas the second order terms 
give rise
to hessians, which are widely
recognised to carry subtle geometric information along the singular locus of 
the theta divisor. For example, in the case of a 
jacobian these hessians
define quadrics containing the canonical image of the associated 
Riemann surface. Or, in the general case one could investigate 
the properties of those principally polarised abelian
varieties that have a singular point of order two on their theta divisor, 
such that the hessian of the theta function at that singular point has 
a certain given rank. This is a recent line of investigation begun by Grushevsky
and Salvati Manni, with interesting connections to the
Schottky problem \cite{gsm} \cite{gsm2}.  

If one moves outside the singular locus of the theta divisor, 
it is not immediately clear whether the
hessian of the theta function continues to have some geometric significance.
In this paper we prove that it does. More
precisely, we show 
that a certain combination of
the gradient and the hessian gives rise to a well-defined theta function 
living on the theta divisor. We can compute its transformation behavior, 
i.e. its order and automorphy factor explicitly. 
As it turns out, the geometric meaning of our function is that 
on the smooth locus of the theta divisor it precisely 
gives the ramification locus of the Gauss map,
i.e. the map sending a smooth point to its tangent space, seen as a point in
the dual of the projectivised tangent space of the origin of the ambient abelian
variety. We can study the situation in more detail for jacobians, where
our theory leads to a global description of
the locus of Weierstrass points on the corresponding Riemann surface in terms of
theta functions and their derivatives. 
By way of example we state and prove an
explicit formula for our theta function in genus two.

\section{Definition and main theorem} \label{mainthm}

% theta, transformation, definition \eta, examples g=1, g=2, main theorem.
% interpretation: support on a single Theta is well-defined, section of a modular
% line bundle.
% mention generic non-vanishing; proof later.

Let $\hh_n$ denote the Siegel upper half space of degree $n>0$. On $\cc^n \times
\hh_n$ we have the Riemann theta function
\[ \theta = \theta(z,\tau) = \sum_{m \in \zz^n} e^{\pi i {}^t m \tau m + 2\pi i
{}^t m z} \, . \] 
Here and henceforth, vectors are column vectors and ${}^t$ denotes transpose.
For any fixed $\tau$, the function $\theta=\theta(z)$ on $\cc^n$ satisfies a
transformation formula
\begin{equation} \label{transformation}
\theta(z+\tau u + v) = e^{ -\pi i {}^t u \tau u - 2 \pi i {}^t u z } \theta(z)
\, ,  
\end{equation}
for all $z$ in $\cc^n$ and for all 
$u,v$ in $\zz^n$. Moreover, it has a symmetry property
\begin{equation} \label{symmetry}
\theta(-z) = \theta(z)  
\end{equation}
for all $z$ in $\cc^n$. Equation (\ref{transformation}) implies that $\div
\theta$ is well-defined as a Cartier divisor on the complex torus
$A=\cc^n/(\zz^n + \tau \zz^n)$. We denote this divisor by
$\Theta$. In fact
$\theta$ gives rise to a global section of $O_A(\Theta)$, and we say that $\theta$
is a ``theta function of order $1$''. By equation
(\ref{symmetry}), the divisor $\Theta$ is symmetric: $[-1]_A^* \Theta = \Theta$.

For $\tau$ varying through $\hh_n$, it turns out that $\theta=\theta(z,\tau)$ is
a modular form of weight $1/2$. More precisely, consider the group 
$\Gamma_{1,2}$ 
of matrices $\gamma = \left( \begin{array}{cc} a & b \\ c & d \end{array} \right)$ in
$\mathrm{Sp}(2n,\zz)$ with $a,b,c,d$ square matrices such that the diagonals of
both ${}^t a c$ and ${}^t b d$ consist of even integers. We can let
$\Gamma_{1,2}$ act
on $\cc^n \times \hh_n$ via
\[ (z,\tau) \mapsto \left( {}^t (c \tau + d)^{-1} z , (a\tau + b)(c \tau
+d)^{-1} \right) \, . \]
Under this action, the Riemann theta function transforms as
\begin{equation} \label{thetatransform}  
\theta( {}^t (c\tau +d)^{-1}z, (a\tau + b)(c \tau + d)^{-1}) =
\zeta_\gamma \det(c\tau +d)^{1/2} 
e^{\pi i {}^t z (c\tau +d)^{-1} c z } \theta(z,\tau)  
\end{equation}
for some $8$-th root of unity $\zeta_\gamma$ (cf. \cite{mum}, p.~189).

We write $\theta_i$ for the first order partial derivative $\partial
\theta/\partial z_i$ and $\theta_{ij}$ for the second order partial derivative
$\partial^2 \theta/\partial z_i \partial z_j$. If $(h_{ij})$ is any square
matrix, we denote by $(h_{ij})^c = (h^c_{ij})$ its cofactor matrix, i.e.
\[ h^c_{ij} = (-1)^{i+j} \det \left( h_{kl} \right)_{k \neq i, l \neq j} \, . \]
The function we want to study in this paper is then the following.
\begin{definition} Let $(\theta_i)$ be the gradient of $\theta$ in the 
$\cc^n$-direction, and
let $(\theta_{ij})$ be its hessian. Then we put
\begin{equation} \label{definition} 
\eta = \eta(z,\tau) = {}^t (\theta_i)   (\theta_{ij})^c   (\theta_j)
\, . 
\end{equation}
We want to consider this as a function on the vanishing locus $\theta^{-1}(0)$
of $\theta$ on $\cc^n \times \hh_n$.
\end{definition}
\begin{example} For $n=1$ we obtain
\[ \eta = \left( \frac{ d \theta}{d z} \right)^2 \, \]
viewed as a function of $(z,\tau)$ in $\cc \times \hh$ with $z \equiv 
(1+\tau)/2 \bmod \zz + \tau \zz$. For $n=2$ we obtain
\[ \eta = \theta_{11} \theta_2^2 - 2 \theta_{12}\theta_1 \theta_2 + \theta_{22}
\theta_1^2 \, , \]
which is already somewhat more complicated.
\end{example}

Our main result is that the function $\eta$ transforms well with
respect to both lattice translations and the action of the congruence symplectic
group $\Gamma_{1,2}$.
\begin{thm} \label{main}
The function $\eta=\eta(z,\tau)$ is a theta function of order $n+1$ and weight
$(n+5)/2$ on the theta divisor. In other words, for any fixed $\tau$ 
in $\hh_n$, the function 
$\eta$ gives rise to a global section of
the line bundle $O_\Theta(\Theta)^{\otimes n+1}$ on $\Theta$ in $ 
A=\cc^n/(\zz^n+\tau \zz^n)$. Furthermore, when viewed as a function of two
variables $(z,\tau)$, the function $\eta$ transforms under the action of
$\Gamma_{1,2}$ with an
automorphy factor $\det( c\tau +d )^{(n+5)/2}$ on $\theta^{-1}(0)$.
\end{thm} 
It follows that for any fixed $\tau$, the zero locus of $\eta$ is well-defined
on $\Theta$. This zero locus contains $\mathrm{Sing} \, \Theta$, 
the singular locus of $\Theta$, as well as the set of $2$-division points $\Theta \cap A[2]$
if $n \geq 2$. Moreover this zero locus 
is stable under the involution $z \mapsto -z$ of
$\Theta$. 

For $\tau$ varying through $\hh_n$ it follows that 
the function $\eta$ gives rise
to a global section of a line bundle
\[ L = O_\Theta(\Theta)^{\otimes n+1} \otimes \pi^* M \]
on the canonical symmetric theta divisor $\Theta$ of a universal principally
polarised abelian variety with level structure $\pi : \mathcal{U}_n \to
\mathcal{A}_n^{(1,2)}= \Gamma_{1,2} \setminus \hh_n$. Here $M$ is a certain line
bundle on $\mathcal{A}_n^{(1,2)}$. It follows from general principles that $M$ is a power of $\lambda$, the
determinant of the Hodge bundle on $\mathcal{A}_n^{(1,2)}$. By counting weights
we find that $M \cong \lambda^{\otimes 2}$.
\begin{example} \label{jacobiderivative}
When $n=1$, the theorem states that 
$ \frac{ d \theta}{d z} $ is of order $1$ and
of modular weight $3/2$ on the theta divisor. 
Both statements follow directly from
(\ref{transformation}) and (\ref{thetatransform}). 
Alternatively, the statement on the modular weight can be seen using Jacobi's
derivative formula (cf. \cite{mum}, p.~64). 
This formula says that
\[ e^{\pi i \tau/4} \frac{ d \theta}{d z} \left( \frac{1+\tau}{2} \right) 
= \pi i \theta \left[ {0 \atop 0} \right](0,\tau)
\theta \left[ {0 \atop 1/2} \right](0,\tau) \theta \left[ {1/2 \atop 0} \right]
(0,\tau)  \]
where $\theta \left[ {0 \atop 0} \right],
\theta \left[ {0 \atop 1/2} \right]$ and $\theta \left[ {1/2 \atop 0} \right]$
are the usual elliptic theta functions with even characteristic. 
Each of the three Thetanullwerte $\theta \left[ {0 \atop 0} \right](0,\tau),
\theta \left[ {0 \atop 1/2} \right](0,\tau)$ and 
$\theta \left[ {1/2 \atop 0} \right](0,\tau)$
is a modular form of weight $1/2$.
\end{example}
A proof of Theorem \ref{main} will be given in Section \ref{proof}.

\section{Properties} \label{properties}

In this section we collect some further properties of $\eta$.
\begin{prop} Assume that $\Theta = \div \theta $ on $ \cc^n/(\zz^n + \tau
\zz^n)$ is non-singular, and denote its canonical bundle by $K_\Theta$. Then
$\eta$ gives rise to a global section $\widetilde{\eta}$ of
$K_\Theta^{\otimes n+1}$, locally given by
\[   \widetilde{\eta}(z) = \eta(z) \cdot \left( (-1)^{i-1}
\frac{dz_1 \wedge \ldots \wedge \widehat{dz_i} \wedge \ldots \wedge
dz_n}{\theta_i(z)}  \right)^{\otimes n+1} \]
wherever $\theta_i(z)$ is non-zero.
\end{prop}
\begin{proof} We view $\eta$ as a global section of
the line bundle $O_\Theta(\Theta)^{\otimes n+1}$ on $\Theta$. 
By the adjunction
formula and the fact that the canonical bundle $K_A$ of $A$ is trivial we have
an identification
\[ O_\Theta (\Theta) \cong (K_A \otimes O_A(\Theta))|_\Theta 
\cong K_\Theta \, .  \]
This identification can be represented locally by the 
Poincar\'e residue map $K_A \otimes O_A(\Theta) \to K_\Theta$ given by
\[ \frac{dz_1 \wedge \ldots \wedge dz_n}{\theta(z)} \longmapsto (-1)^{i-1}
\frac{dz_1 \wedge \ldots \wedge \widehat{dz_i} \wedge \ldots \wedge
dz_n}{\theta_i(z)}  \]
wherever $\theta_i(z) \neq 0$. From this the corollary follows.
\end{proof}
By way of illustration, we 
make the multidifferential $\widetilde{\eta}$ more explicit in the
case that $n=2$. We need the notions of a Wronskian
differential and of Weierstrass points, which we recall briefly. Let $X$ be a
compact Riemann surface of genus $n>0$, and let $\zeta=(\zeta_1,\ldots,\zeta_n)$
be any basis of $H^0(X,K_X)$. The Wronskian differential $\omega_\zeta$ of $\zeta$
is then given as follows: let $t$ be a local coordinate, and write $\omega_j(t)
=f_j(t) dt$ with the $f_j$ holomorphic. We put
\[ \omega_\zeta (t) = \det \left( \frac{1}{(i-1)!} \frac{d^{i-1} f_j(t)}{dt^{i-1}}
\right) \cdot (dt)^{\otimes n(n+1)/2 } \, . \]
This local definition gives, in fact, rise to a global section of the
line bundle $K_X^{\otimes n(n+1)/2}$. It can be proved that this section is
non-zero. If the basis $\zeta$ is changed, the differential $\omega_\zeta$
changes by a non-zero scalar. Hence the divisor $W=\div \omega_\zeta$ is
independent of the choice of $\zeta$. We call this divisor the
classical divisor of Weierstrass points of $X$. It has degree $n^3-n$.   
\begin{example} \label{wronskiangenustwo}
Assume that $A=\cc^2/(\zz^2 + \tau \zz^2)$ is an indecomposable
abelian surface with 
theta divisor $\Theta = \div \theta$.
Then $\Theta$ is a compact Riemann surface of genus $2$, and there is a canonical
identification $H^0(A,\Omega^1) \cong H^0(\Theta,K_\Theta)$. Let
$\zeta=(\zeta_1(z),\zeta_2(z))$ be the basis of $H^0(\Theta,K_\Theta)$ 
corresponding under this identification to
the standard basis $(dz_1,dz_2)$ of $H^0(A,\Omega^1)$. We claim that
$\widetilde{\eta}(z)$ is equal to the Wronskian differential $\omega_\zeta$ of
$\zeta$. This amounts to a small computation: choose an open subset of $\Theta$
where $z_1$ is a local coordinate. In this local coordinate we can write 
$\zeta_1(z) = dz_1$, $\zeta_2(z) =
z'_2(z) \cdot dz_1$, with $'$ denoting derivative with respect to $z_1$, so that
\[ \omega_\zeta(z) = \det \left( \begin{array}{cc}
	1  &  z'_2(z) \\
	0  &	z''_2(z) \end{array} \right) 
	\cdot (dz_1)^{\otimes 3} = 
	z''_2(z) \cdot (dz_1)^{\otimes 3} \, . \]
Now, since for $z$ on $\Theta$ we have
\[ \theta_1(z) dz_1 + \theta_2(z) dz_2 = 0 \, , \]
the formula 
\[ z'_2(z) = - \frac{ \theta_1(z)}{\theta_2(z)} \] 
holds, leading to
\begin{eqnarray*} 
 z''_2(z) & = & - \frac{\theta'_1(z) \theta_2(z) -
 \theta'_2(z) \theta_1(z) }{ \theta_2(z)^2} \\
   & = & - \frac{\left( \theta_{11}(z) + 
   \theta_{12}(z) \cdot - \frac{\theta_1(z)}{\theta_2(z)} \right) \theta_2(z) -
   \left( \theta_{12}(z) + \theta_{22}(z) 
   \cdot - \frac{\theta_1(z)}{\theta_2(z)} \right) \theta_1(z) }{ \theta_2(z)^2
   }  \\
   & = & - \frac{ \theta_{11} (z) \theta_2^2(z) - 2 \theta_{12}(z) \theta_1(z)
\theta_2(z) +
\theta_{22}(z) \theta_1^2(z)  }{ \theta_2(z)^3     } \\
   & = & - \frac{ \eta(z) }{ \theta_2(z)^3 } \, . 
\end{eqnarray*}
We find, indeed, 
\[ \omega_\zeta(z) = - \frac{ \eta(z) }{ \theta_2(z)^3 } \cdot (dz_1)^{\otimes 3} =
\widetilde{\eta}(z) \,
. \]
A similar computation can be done on the locus where $z_2$ is a local
coordinate, giving $\omega_\zeta(z) = \widetilde{\eta}(z)$ 
globally on $\Theta$, as required. 
Note that our identification $\omega_\zeta(z) = \widetilde{\eta}(z)$ gives, 
as a corollary, that $\div \eta=W$, the divisor of
Weierstrass points of $\Theta$.
In Section \ref{jacobians} we will prove a generalisation of this result.
\end{example}
It is interesting to know when $\eta$ is identically zero on the theta
divisor.
\begin{prop} \label{decomposables}
If $A$ is a decomposable principally polarised abelian
variety, then $\eta$ is identically zero on the theta divisor.
\end{prop} 
\begin{proof} Let us suppose that 
$A = A_1 \times A_2$ with $A_1$ given by a matrix $\tau_1$ in 
$\hh_k$ and $A_2$ given by a matrix $\tau_2$ in $\hh_{n-k}$ where $k$ is an
integer 
with $0<k<n$.
We can write
$\theta(z)=F(z_1,\ldots,z_k)G(z_{k+1},\ldots,z_n)$ where $F,G$ are 
the Riemann theta functions for $A_1$ and $A_2$, respectively. 
Let $\Theta_1 \subset A_1$ be the divisor of $F$, and $\Theta_2 \subset A_2$ 
the divisor of $G$. By symmetry, it suffices to prove that $\eta$ is zero on
$\Theta_1 \times A_2 \subset \Theta = \div \theta$. 
On $\cc^n$ we have
\[ {}^t (\theta_i) = (F_iG \, , \, FG_i)  \]
and
\[  (\theta_{ij}) = \left( \begin{array}{cc} 
F_{ij}G &  F_iG_j \\
F_jG_i & FG_{ij} \end{array} \right) \, . \]            
The subset $\Theta_1 \times A_2 \subset \Theta$ is given by the vanishing of
$F$; there we find
\[ {}^t (\theta_i) = (F_iG , 0 ) \quad , \quad 
(\theta_{ij}) = \left( \begin{array}{cc} 
F_{ij}G &  F_iG_j \\
F_jG_i & 0 \end{array} \right) \, . \]
Note that both $F_iG_j$ and $F_jG_i$, being a product of a vector and a
covector, have rank $\leq 1$. Hence, a minor at
$(i,j)$ with $1 \leq i,j \leq k$ in $(\theta_{ij})$ has rank $\leq 1+(k-1) 
=k$. If $k < n-1$, the determinant of this minor vanishes, and 
the cofactor matrix of $(\theta_{ij})$ has the shape
\[ (\theta_{ij})^c = \left( \begin{array}{cc} 0 & * \\ * & * \end{array} 
\right) \, .
\]
We find then that
\[ \eta = {}^t (\theta_i)   (\theta_{ij})^c  (\theta_j) =
( * , 0 )   \left( \begin{array}{cc} 0 & * \\ * & * \end{array} \right)
  \left(
\begin{array}{c} * \\ 0 \end{array} \right) = 0 \, . \]
If $k=n-1$, the last row of $(\theta_{ij})$ and the vector 
${}^t (\theta_i)$ are linearly
dependent. If $1 \leq i \leq n-1$ we see that the $i$-th entry of
$(\theta_{ij})^c (\theta_j)$ is just the determinant of 
the matrix obtained from
$(\theta_{ij})$ by removing the $i$-th row and adding ${}^t (\theta_i)$ in its
place. But this matrix then contains two linearly dependent vectors, and its
determinant vanishes. So we obtain
\[ \eta =  {}^t (\theta_i)   (\theta_{ij})^c  (\theta_j) = (*,0) \left(
\begin{array}{c} 0 \\ * \end{array} \right) = 0 \]
in this case as well.
\end{proof}

\section{Interpretation} \label{interpretation}

The contents of the present section are based on kind suggestions made by 
Prof. Ciro Ciliberto.

As we observed in Section \ref{mainthm}, for any
complex principally polarised abelian variety the zero locus of $\eta$ is
well-defined on the theta divisor. Given the simple
description of $\eta$ one expects that this zero locus has some intrinsic
geometric interpretation. This indeed turns out to be the case.
\begin{thm} \label{ciro}
Let $(A,\Theta)$ be a complex principally polarised 
abelian variety. On the
smooth locus $\Theta^s$ of $\Theta$, the zero locus 
of $\eta$ is precisely 
the ramification locus of the Gauss map
\[ \Gamma : \Theta^s \longrightarrow \pp( T_0 A)^*  \]
sending a point on $\Theta^s$ to its tangent space, translated to a subspace of
$T_0 A$.
\end{thm} 
\begin{proof} It follows from formula (\ref{definition}) that a point $x$ on
$\Theta^s$ is in the zero locus of $\eta$ if and only if the quadric $Q$ in
$\pp(T_0 A)$ defined by the hessian is tangent to the projectivised tangent 
hyperplane $\pp(T_x \Theta)$
defined by the gradient. The latter condition is equivalent to the condition 
that $Q$ when restricted to $\pp(T_x \Theta)$ becomes degenerate. 
Now note that $Q$ when viewed as a linear map $Q :
T_x \Theta \to (T_x \Theta)^*$ can be identified with the tangent
map $d \Gamma : T_x \Theta \to T_{\Gamma(x)}(\pp(T_0A)^*) =
(T_x \Theta)^*$ of $\Gamma$. 
The locus where this 
map is degenerate is precisely the ramification locus of $\Gamma$. 
\end{proof}
Using the above interpretation, a converse to Proposition \ref{decomposables} can be
readily proved.
\begin{cor} \label{finaldecomposables} We have that $\eta$ is identically zero on the
theta divisor if and only if $A$ is a decomposable abelian variety.
\end{cor} 
\begin{proof} We need to prove that if $\Theta$ is irreducible, then 
the Gauss map on $\Theta$ has a proper
ramification locus. But according to \cite{kempf}, Corollary 9.11
the Gauss map is generically finite and dominant 
in this case, and the result follows.
\end{proof}
\begin{rem} It follows from Theorem \ref{main} that in the indecomposable case 
the divisor of $\eta$ belongs to the linear
system defined by $(n+1)\Theta$ on $\Theta$. This fact
can also be explained as follows. For
simplicity, let us assume that $\Theta$ is non-singular. The Hurwitz formula
applied to the Gauss map $\Gamma : \Theta \to  \pp=\pp( T_0 A)^* $ gives that
$K_\Theta = \Gamma^* K_\pp + R$, where $R$ is the ramification locus of
$\Gamma$. By definition $\Gamma^*
O_\pp(1) \cong K_\Theta$ and hence $\Gamma^*K_\pp \cong K_\Theta^{\otimes -n}$.
We find that $R$ is in the linear system belonging to 
$K_\Theta^{\otimes n+1}$. As
we have seen at the end of the previous section, this is the same as
$O_\Theta(\Theta)^{\otimes n+1}$. Alternatively, this remark shows that the
equations $\theta = \eta= 0$ define the right scheme structure on $R$.
\end{rem}
\begin{rem} The above description of the vanishing locus of $\eta$ can also be
found in the paper \cite{gsm2} of Grushevsky and Salvati Manni. In that paper
they give an application of the form $\eta$ to the study of certain
codimension-$2$ cycles in the moduli space of principally polarised abelian
varieties. In particular they are able to give a moduli interpretation of a
certain cycle $\mathrm{R}_g$ introduced by Debarre in \cite{deb}, Section~4. 
\end{rem}

\section{Proof of the main theorem} \label{proof}

% the proof

In this section we give a proof of Theorem \ref{main}. We start with the
statement on the order of $\eta$. We fix an element $\tau$ in $\hh_n$. The
case $n=1$ being discussed already 
in Example \ref{jacobiderivative}, we assume here that $n \geq 2$.
Recall from equation (\ref{transformation}) that
\begin{equation} \label{transformation2}
\theta(z+\tau u + v) = p(z,u) \theta(z)
\end{equation}
for all $z$ in $\cc^n$ and all $u,v$ in $\zz^n$, where
\[ p(z,u) = e^{ -\pi i {}^t u \tau u - 2 \pi i {}^t u z } \, . \]
We need to prove that
\[ \eta(z+\tau u + v) = p(z,u)^{n+1} \eta(z) \]
for all $z$ in $\cc^n$ with $\theta(z)=0$. Denote by $p_i$ the first order
partial derivative $\partial p /\partial z_i$. From (\ref{transformation2}) we
have, for $z$ with $\theta(z)=0$,
\[ \theta_i(z+\tau u + v) = p(z,u) \theta_i(z) \, \]
and
\[ \theta_{ij} (z+ \tau u + v) = p(z,u) \theta_{ij}(z) + p_i(z,u) \theta_j(z) +
p_j(z,u) \theta_i(z) \, . \]
We are done, therefore, if we can prove, formally in some domain $R$
containing the symbols $p,p_i,\theta_i,\theta_{ij}$, the identity
\[ {}^t (p \theta_i) (p \theta_{ij} + p_i \theta_j + p_j \theta_i)^c (p
\theta_j) = p^{n+1} {}^t (\theta_i) (\theta_{ij})^c (\theta_j) \, \]
or, equivalently, the identity
\[ {}^t (\theta_i) (p \theta_{ij} + p_i \theta_j + p_j \theta_i)^c (\theta_j) =
{}^t (\theta_i) (p \theta_{ij})^c (\theta_j) \, . \]
At this point we introduce some notation. Let $h=(h_{i_kj_l})$ be an 
$R$-valued matrix with rows and columns indexed by
finite length ordered integer tuples $I=(\ldots,i_k,\ldots)$ and $J=(\ldots,j_l,\ldots)$,
respectively. 
If $I' \subset I$ and $J' \subset J$ are proper
subtuples, we denote by $h^{I'}_{J'}$ the submatrix obtained 
from $h$ by deleting the
rows indexed by $I'$ and the columns indexed by $J'$. If $h$ is a square
$R$-valued matrix with both the row index set
$I=(i_1,\ldots,i_m)$ and the column index set 
$J=(j_1,\ldots,j_m)$ subtuples of $(1,2,\ldots,n)$,
we define $\eta(h)$ to be the
element
\[ \eta(h) = {}^t (\theta_{i_k}) (h_{i_k j_l})^c (\theta_{j_l})  \]
of $R$. It can be written more elaborately as
\[ \eta(h) = \sum_{k=1}^m \sum_{l=1}^m (-1)^{k+l} \theta_{i_k}
\theta_{j_l}  
\det h^{\{i_k\}}_{\{j_l\}} \, . \]  
On the other hand, 
the identity we need to prove can be written more compactly as
\[ \eta(p \theta_{ij} + p_i \theta_j + p_j \theta_i  ) = \eta(p \theta_{ij} ) \,
. \]
The function $\eta$ can be defined recursively.
\begin{lem} \label{identity}
Let $h=(h_{i_kj_l})$ be a square $R$-valued matrix of size $m \geq 2$ with 
both the row index set
$I=(i_1,\ldots,i_m)$ and the column index set 
$J=(j_1,\ldots,j_m)$ subtuples of $(1,2,\ldots,n)$. 
Then the identity
\[ (m-1) \eta(h) = \sum_{k=1}^m \sum_{l=1}^m (-1)^{k+l} h_{i_k j_l}
\eta\left( h^{ \{
i_k \} }_{ \{ j_l \} } \right) \]
holds.
\end{lem}
\begin{proof} The double sum on the right hand side can be expanded as
\begin{eqnarray*} 
\sum_{(k,l)} (-1)^{k+l} h_{i_k j_l}
\eta\left( h^{ \{
i_k \} }_{ \{ j_l \} } \right) & = & 
\sum_{(k,l)} (-1)^{k+l} h_{i_kj_l} \sum_{(k',l') \neq (k,l)} (-1)^{k'+l'}
\theta_{i_{k'}} \theta_{j_{l'}} \det h^{ \{
i_k , i_{k'} \} }_{ \{ j_l , j_{l'} \} } \\
 &= & \sum_{(k',l')} (-1)^{k'+l'} \theta_{i_{k'}} \theta_{j_{l'}} \sum_{(k,l) \neq
 (k',l')} (-1)^{k+l} h_{i_k j_l} \det h^{ \{
i_{k'} , i_k \} }_{ \{ j_{l'} , j_l \} } \, .
\end{eqnarray*}
Expanding the determinant of $h^{ \{ i_{k'} \} }_{ \{ j_{l'} \} } $ 
along each of its
rows, one finds
\[ (m-1) \det h^{ \{ i_{k'} \} }_{ \{ j_{l'} \} } = 
\sum_{(k,l) \neq (k',l')} (-1)^{k+l}
h_{i_k j_l} \det h^{ \{
i_{k'} , i_k \} }_{ \{ j_{l'} , j_l \} } \, . \]
Combining both formulas one gets
\[ \sum_{(k,l)} (-1)^{k+l} h_{i_k j_l}
\eta\left( h^{ \{
i_k \} }_{ \{ j_l \} } \right) = (m-1) \sum_{(k',l')} (-1)^{k'+l'}
\theta_{i_{k'}}
\theta_{j_{l'}} \det h^{ \{ i_{k'} \} }_{ \{ j_{l'} \} } = (m-1) \eta(h) \, , \]
as required.
\end{proof}
Denote by $a=(a_{ij})$ the $n$-by-$n$ $R$-valued matrix 
\[ (a_{ij}) = ( p_i \theta_j + p_j \theta_i ) \, . \]
It has the property that the $\eta$'s of all its square submatrices
vanish.
\begin{lem} \label{vanishing_of_eta}
Let $h$ be any square submatrix of $a$ of size $\geq 2$. Then
$\eta(h)=0$.
\end{lem} 
\begin{proof} By Lemma \ref{identity}, the statement follows by induction once
we prove the special case that $h$ is a square submatrix of size $2$. In this
case $h$ has the shape 
\[ \left( \begin{array}{cc}
	p_{i_1} \theta_{j_1} + p_{j_1} \theta_{i_1}  & 
		p_{i_1} \theta_{j_2} + p_{j_2} \theta_{i_1} \\
	p_{i_2} \theta_{j_1} + p_{j_1} \theta_{i_2}  &
		p_{i_2} \theta_{j_2} + p_{j_2} \theta_{i_2}
	\end{array} \right) \, , \]
and $\eta(h)$ has an expansion
\begin{eqnarray*} 
\eta(h) & =  \theta_{i_1}\theta_{j_1} 
	\left( p_{i_2} \theta_{j_2} + p_{j_2} \theta_{i_2} \right) -
	\theta_{i_1}\theta_{j_2} 
	\left( p_{i_2} \theta_{j_1} + p_{j_1} \theta_{i_2} \right) \\
	& - \theta_{i_2}\theta_{j_1} 
	\left( p_{i_1} \theta_{j_2} + p_{j_2} \theta_{i_1} \right) +
	\theta_{i_2}\theta_{j_2} 
	\left( p_{i_1} \theta_{j_1} + p_{j_1} \theta_{i_1} \right) \, .
\end{eqnarray*} 
This is identically equal to zero.
\end{proof}
Denote by $b=(b_{ij})$ the $n$-by-$n$ $R$-valued matrix 
\[ (b_{ij}) = (p \theta_{ij} ) \, , \]
and by $\widetilde{b}=(\widetilde{b}_{ij})$ the $n$-by-$n$ $R$-valued matrix 
\[ (\widetilde{b}_{ij}) = ( p \theta_{ij} + p_i \theta_j + p_j \theta_i ) \, .
\]
We need to prove that $\eta(\widetilde{b}) = \eta(b)$. We can expand 
$\eta(\widetilde{b}) - \eta(b) $ 
as
\begin{equation} \label{expandetab} 
\eta(\widetilde{b}) - \eta(b) = 
\sum_{m=0}^{n-2} p^m \sum_{ (I,J)  } \varepsilon_{I,J}
\theta_{i_1 j_1} \cdots \theta_{i_m j_m} \eta(a^I_J) \, . 
\end{equation}
Here the second sum is over all pairs $(I,J)$ of subtuples of $(1,2\ldots,n)$ of
length $m$, and $\varepsilon_{I,J}$ is a sign. 
In order to see this, expand 
in the cofactor matrix $\widetilde{b}^c $ of $\widetilde{b}$ each minor 
as a sum of $(n-1)!$ terms. A product 
$p^m \theta_{i_1 j_1} \cdots \theta_{i_m j_m} $  
occurs as a factor in such a
term exactly at all entries $(k,l)$ of $\widetilde{b}^c $ for which $(k,l)$ is
not in $I \times J$. If $(k,l)$ is such an entry, at that entry the product 
$p^m \theta_{i_1 j_1} \cdots \theta_{i_m j_m} $ is multiplied, up to a sign
$\varepsilon_{I,J}$ depending only on $I,J$, by 
$\det \left( a^{I \cup \{k \} }_{J \cup \{ l \} } \right)$. 
This determinant is understood to be equal to
$1$ if $m=n-1$. It follows that
the entry $(k,l)$ contributes to $\eta(\widetilde{b})$ with a term
\[ p^m \theta_{i_1 j_1} \cdots \theta_{i_m j_m} \varepsilon_{I,J} 
(-1)^{k+l} \theta_k \theta_l  
\det \left( a^{I \cup \{k \} }_{J \cup \{ l \} } \right) \, . \] 
Summing over all possible $(k,l)$
we obtain, if $m<n-1$, a contribution 
$ p^m \theta_{i_1 j_1} \cdots \theta_{i_m j_m} \varepsilon_{I,J} \eta(a^I_J)  $,
and if $m=n-1$ the contribution $\eta(b)$. 

By Lemma \ref{vanishing_of_eta}, every $\eta(a^I_J)$ with $I,J$ of size smaller
than $n-1$ is zero. Therefore 
all terms in the summation on the right hand side in (\ref{expandetab}) 
vanish. This proves the first half of Theorem \ref{main}.

In order to prove the statement on the modular weight, we recall from equation
(\ref{thetatransform}) that 
\begin{equation} \label{thetatransformbis} 
 \theta( {}^t (c\tau +d)^{-1}z, (a\tau + b)(c \tau + d)^{-1}) =
\zeta_\gamma \det(c\tau +d)^{1/2} q(z,\gamma,\tau)
 \theta(z,\tau) 
\end{equation} 
for all $z$ in $\hh_n$ and all 
$\gamma = \left( \begin{array}{cc} a & b \\ c & d \end{array}
\right)$ in $\Gamma_{1,2}$, where
\[ q(z,\gamma,\tau) = e^{\pi i {}^t z (c\tau +d)^{-1} c z } \]
and where $\zeta_\gamma$ is an $8$-th root of unity. We claim that 
\begin{equation} \label{etatransform} 
\eta( {}^t (c\tau +d)^{-1}z, (a\tau + b)(c \tau + d)^{-1}) =
\det(c\tau +d)^{(n+5)/2} \zeta_\gamma^{n+1} q(z,\gamma,\tau)^{n+1} \eta(z,\tau) 
\end{equation}
for all $(z,\tau)$ satisfying $\theta(z,\tau)=0$. This is just a calculation. 
For $(z,\tau)$ with
$\theta(z,\tau)=0$ we have by (\ref{thetatransformbis})
\[ \left( 
\theta_i({}^t (c\tau +d)^{-1}z, (a\tau + b)(c \tau + d)^{-1}) \right) = 
{}^t (c\tau
+d) \zeta_\gamma \det(c\tau +d)^{1/2} q(z,\gamma,\tau) 
\left( \theta_i(z,\tau) \right)  \, . \]
Furthermore, for such $(z,\tau)$ we have
\begin{eqnarray*} & \left( \theta_{ij}({}^t (c\tau +d)^{-1}z, 
(a\tau + b)(c \tau + d)^{-1}) \right) = \hspace{5cm} \\
& {}^t (c\tau
+d) \zeta_\gamma \det(c\tau +d)^{1/2} \left( q(z,\gamma,\tau) 
\theta_{ij}(z,\tau) +
q_i(z,\gamma,\tau)
\theta_j(z,\tau) + q_j(z,\gamma,\tau)\theta_i(z,\tau)  \right) (c\tau + d) \, .
\end{eqnarray*}
We can write, at least for the purpose of this proof,
\begin{equation} \label{alternative}
\eta = \det (\theta_{ij}) {}^t (\theta_i) (\theta_{ij})^{-1} (\theta_j) \, . 
\end{equation}
This gives that
$ \eta( {}^t (c\tau +d)^{-1}z, (a\tau + b)(c \tau + d)^{-1})$ is equal to
\begin{eqnarray*} & & \det \left( {}^t (c\tau +d) \zeta_\gamma 
\det(c\tau +d)^{1/2} \left( q  \theta_{ij}(z,\tau) +
q_i \theta_j(z,\tau) + q_j \theta_i(z,\tau)  \right) (c\tau + d)
\right)  \cdot \\
& & \det(c\tau +d)^{1/2} \, \zeta_\gamma {}^t \left( q  \theta_i(z,\tau) \right) 
(c\tau +d) \cdot \\
& &(c\tau + d)^{-1} \left( q  \theta_{ij}(z,\tau) +
q_i \theta_j(z,\tau) + q_j \theta_i(z,\tau)  \right)^{-1}  
\zeta_\gamma^{-1} \det(c\tau +d)^{-1/2} \, {}^t (c\tau +d)^{-1}  \cdot  \\
& & {}^t (c\tau+d) \zeta_\gamma \det(c\tau +d)^{1/2} \,  
\left( q  \theta_j(z,\tau) \right) \, .  
\end{eqnarray*}
This simplifies to 
\begin{eqnarray*} 
&   \det(c\tau +d)^{2+\frac{n}{2}+\frac{1}{2}} \, \zeta_\gamma^{n+1} \cdot 
\hspace{3cm}
\\
&  {}^t(q \theta_i(z,\tau)) \left( q  \theta_{ij}(z,\tau) +
q_i \theta_j(z,\tau) + q_j \theta_i(z,\tau)    \right)^c
(q  \theta_j(z,\tau) ) \, , 
\end{eqnarray*}
which, in turn, is equal to
\[ \det(c\tau + d)^{(n+5)/2} 
\zeta_\gamma^{n+1} q^{n+1} \eta(z,\tau) \]
by the same methods as above when we dealt with the order of $\eta$. This
completes the proof of Theorem \ref{main}.

\section{Jacobians} \label{jacobians}

% non-vanishing for jacobians, in fact: Weierstrass points are given
% mention at the end once more the problem of finding a good geometric
% interpretation in general
% also, refer back to the determination of M in one of the corollaries

The purpose of this section is to study $\eta$ for jacobians. 

Assume that
$(A,\Theta)=(\cc^n/(\zz^n + \tau \zz^n),\div \theta)$ with $n \geq 2$ is the
jacobian belonging to a compact Riemann surface $X$ marked with a symplectic
basis $\mathcal{B}=(A_1,\ldots,A_n,B_1,\ldots,B_n)$ of homology. This situation
determines uniquely a basis $(\zeta_1,\ldots,\zeta_n)=\zeta= \zeta_{\mathcal{B}}$
of $H^0(X,K_X)$ such that
\[ \int_{A_i} \zeta_j = \delta_{ij} \quad , \quad \int_{B_i} \zeta_j = \tau_{ij}
\, . \]
The isomorphism of $\cc$-vector spaces $H^0(X,K_X)^* \isom \cc^n$ given by
$\zeta$ gives an isomorphism of complex tori
\[ H^0(X,K_X)^*/H_1(X,\zz) \isom \cc^n/(\zz^n + \tau \zz^n) = A \, . \]
By a theorem of Abel-Jacobi, the natural map
\[ AJ : \mathrm{Pic}^0 X \longrightarrow \cc^n/(\zz^n + \tau \zz^n) \quad ,
\quad \sum (P_i - Q_i) \mapsto \sum \int_{Q_i}^{P_i} {}^t
(\zeta_1,\ldots,\zeta_n) \]
is a bijection. By a theorem of Riemann, there is a unique element
$\Delta=\Delta_\mathcal{B}$ of $\mathrm{Pic}^{1-n} X $ such that under the
composition of bijections
\[ \mathrm{Pic}^{n-1} X \, {\buildrel t_\Delta \over \longrightarrow} \, 
\mathrm{Pic}^0 X \, 
{\buildrel AJ \over \longrightarrow} \, \cc^n/(\zz^n + \tau \zz^n) = A \, , \]
the set $\Theta_0 = \{ [D] \in \mathrm{Pic}^{n-1} X \, : \, h^0(D) > 0 \}$ is
identified with $\Theta = \div \theta$ on $A$. From now on, we will take this
identification of $(A,\Theta)$ with $(\mathrm{Pic}^{n-1} X, \Theta_0)$ for
granted. 

We have a natural surjection 
$\Sigma : X^{(n-1)} \to \Theta$ which is an isomorphism
above $\Theta^s$. If we let this isomorphism be
followed by the Gauss map $\Gamma : \Theta^s \to \pp(T_0 A)^*=\pp(H^0(X,K_X))$
we get the map sending a non-special 
divisor $D$ of degree $n-1$ on $X$ to the linear span of its points on the
canonical image in $\pp(H^0(X,K_X)^*)$. By Riemann-Roch, this span is indeed a
hyperplane. 

We denote by
\[ \kappa : X \longrightarrow \Theta \subset A \]
the map sending $x$ to the class of $(n-1) \cdot x$. 
It is natural to study the pullback of $\eta$ along $\kappa$. We claim the
following result:
\begin{thm} \label{jaco} The section $\kappa^*\eta$ is not identically zero. The
divisor of $\kappa^*\eta$ is equal to $(n-1)W$, where $W$ is the divisor of
Weierstrass points of $X$.
\end{thm}
We proceed in a few steps, starting with some notation. 
For any $(m,i)$ in $\zz \times \zz_{\geq 0}$ we put
\[ f_{m,i}(x) = \left. \left( \frac{d^i}{dy^i} \theta( (n-1-m)x + my) \right) 
\right|_{y=x} \, , \]
interpreted as a global section of the bundle $\kappa^* O_A(\Theta) \otimes
K_X^{\otimes i}$ of differential $i$-forms with coefficients in $\kappa^*
O_A(\Theta)$. For example, for $i=0$ we get
\[ f_{m,0}(x) = 0 \]
identically, for $i=1$ we get
\[ f_{m,1}(x) = m \sum_{j=1}^n 
\theta_j(\kappa(x))
\zeta_j (x) \, , \]
and for $i=2$:
\[ f_{m,2}(x) = m^2 \sum_{j,k=1}^n \theta_{jk}(\kappa(x))\zeta_j(x) \zeta_k(x) +
m \sum_{j=1}^n \theta_j(\kappa(x)) \zeta'_j(x) \, . 
\]
For $m=0,\ldots,n-1$ all $f_{m,i}$ are identically equal to zero. Of particular
interest for us will be the section
\[ F(x) = \frac{1}{n!} f_{-1,n}(x)   \]
of $\kappa^*O_A(\Theta) \otimes K_X^{\otimes n}$. It turns out that $\div F =
W$.
\begin{lem} \label{fayslemma} 
(Cf. \cite{jor}, Corollary 3)
The section $F$ of $\kappa^*O_A(\Theta) \otimes K^{\otimes n}_X$ 
is not identically zero, and we have $\div F = W$. 
\end{lem}
\begin{proof} Consider the map 
$\Phi : X \times X \to \mathrm{Pic}^{n-1} X$ given by 
$(x,y) \mapsto nx-y$.
From \cite{fay},  p.~31 we obtain that $\Phi^* \Theta = W \times X + n \cdot
\Delta_X$. Restricting to the diagonal we get
\[ \kappa^* \Theta = \Phi^* \Theta|_{\Delta_X} = \left( W \times X + n \cdot
\Delta_X \right) |_{\Delta_X} \]
and by the adjunction formula it follows that
\[ \kappa^*O_A(\Theta) \isom O_X(W) \otimes K_X^{\otimes -n} \]
via $F(x) \cdot (dx)^{\otimes - n} \mapsto 1_W \otimes  (dx)^{\otimes - n}$.
Here $1_W$ denotes the tautological section of $O_X(W)$.
\end{proof} 
\begin{lem} \label{ramif}
We have $\kappa^*\eta(x)=0$ if and only if $x$ is a Weierstrass
point. In particular $\kappa^*\eta$ is not identically zero.
\end{lem}
\begin{proof} According to \cite{deb}, p.~691 the ramification locus of the
Gauss map is precisely given by the set of divisors $D+x$ with $D$ effective of
degree $g-2$ and $x$ a point of $X$ such that $D+2x$ is dominated by a canonical
divisor, i.e. such that $K_X-D-2x$ is linearly equivalent to an effective
divisor. Thus, to say that $\eta(\kappa(x))=0$ means precisely that $h^0(K_X-n
\cdot x) >0$ or equivalently, by Riemann-Roch, that $h^0(n \cdot x) >1$. But
this means precisely that $x$ is a Weierstrass point.
\end{proof} 
\begin{proof}[Proof of Theorem \ref{jaco}] It follows from Section \ref{mainthm}
that $\kappa^*\eta$ is a global section of $\kappa^*O_A(\Theta)^{\otimes n+1}
\otimes \lambda^{\otimes 2}$ where $\lambda$ is the trivial bundle $\det
H^0(X,K_X) \otimes O_X$. By Lemma \ref{ramif} we know that $\kappa^*\eta$ is
non-zero. It is stated in Lemma \ref{fayslemma} that $F$ is a non-zero global
section of $\kappa^*O_A(\Theta) \otimes K^{\otimes n}_X$. As was observed by
Arakelov (cf. \cite{ar}, Lemma 3.3), the bundle
\[ K_X^{\otimes n(n+1)/2} \otimes \lambda^{\otimes -1} \] has a non-zero section
given by 
\[   \xi_1 \wedge \ldots \wedge \xi_n \mapsto \frac{\xi_1 \wedge \ldots
\wedge \xi_n}{\zeta_1 \wedge \ldots \wedge \zeta_n} \cdot \omega_\zeta   \]
with $\omega_\zeta$ the Wronskian differential on $\zeta$. Combining, we find
that $\omega_\zeta^2 \otimes \kappa^*\eta \otimes F^{\otimes -(n+1)}$ is a
non-zero global section of $O_X$. Hence it is a non-zero constant. We find
\[ \div \kappa^*\eta = (n+1) \div F - 2 \div \omega_\zeta = (n+1)W-2W=(n-1)W \]
as required.
\end{proof}
\begin{rem} \label{formula} An elaborate computation shows that actually
\[  \omega_\zeta^2 \otimes \kappa^*\eta = F^{\otimes n+1} \, . \]
Let us prove this relation in the case that $n=2$. So we look at 
indecomposable $(A,\Theta)$ with 
$A=\cc^2/(\zz^2+\tau \zz^2)$ and $\Theta =X= \div \theta$. 
For $z$ on $\Theta =X$ put
\[ P = \theta_{11}(z) \zeta_1(z)^{\otimes 2} 
+ 2 \theta_{12}(z)\zeta_1(z) \otimes \zeta_2(z)
+ \theta_{22}(z) \zeta_2(z)^{\otimes 2} \]
and
\[ Q = \theta_1(z) \zeta'_1(z) + \theta_2(z) \zeta'_2(z) \, . \]
Then according to what we have said before Lemma \ref{fayslemma} we have
\[ P+Q=0 \quad , \quad P-Q = 2F \, . \]
We conclude that
\[ F = P =  \theta_{11}(z) \zeta_1(z)^{\otimes 2} 
+ 2 \theta_{12}(z)\zeta_1(z) \otimes \zeta_2(z)
+ \theta_{22}(z) \zeta_2(z)^{\otimes 2} \]
and writing as before
\[ \zeta_1(z) = dz_1 \quad , \quad \zeta_2(z) = dz_2 = z_2'(z) dz_1 = - 
\frac{\theta_1(z)}{\theta_2(z)} dz_1 \] 
we get
\begin{eqnarray*} 
F(z) & = &  \theta_{11}(z) \zeta_1(z)^{\otimes 2} 
+ 2 \theta_{12}(z)\zeta_1(z) \otimes \zeta_2(z)
+ \theta_{22}(z) \zeta_2(z)^{\otimes 2} \\
  & = & \left( 
  \theta_{11}(z) + 2 \theta_{12}(z) z'_2(z) + \theta_{22}(z) (z'_2(z))^2 \right)
  \cdot (dz_1)^{\otimes 2} \\
  & = &   \frac{ \theta_{11}(z) \theta_2(z)^2 - 2 \theta_{12} \theta_1(z)
  \theta_2(z) + \theta_{22}(z) \theta_1(z)^2 }{ \theta_2(z)^2 }  
  \cdot (dz_1)^{\otimes 2} \\
  & = & \frac{ \eta(z) }{ \theta_2(z)^2 } \cdot (dz_1)^{\otimes 2} \, . 
\end{eqnarray*}
We have seen in Example \ref{wronskiangenustwo}
that $\omega_\zeta(z) = - (\eta(z)/ \theta_2(z)^3) (dz_1)^{\otimes
3}$. Combining we find
\[  \omega_\zeta^2 \otimes \kappa^*\eta = F^{\otimes 3}  \]
as required. We note in passing that our formula (for general $n$) leads to an
alternative description of one of the analytic invariants studied in \cite{jong}. 
\end{rem}

\section{Explicit formula} \label{explicit}

In the case that $n=2$ it is possible to give 
a closed formula for $\eta$ using a more familiar theta function. 
This formula can
be viewed as a generalisation of Jacobi's derivative formula (cf. Example
\ref{jacobiderivative}), which gives
$\eta$ in the case that $n=1$ as a product of even Thetanullwerte. Recall that for
$a,b$ column vectors of dimension $n$ 
with entries in $\{0,1/2\}$ we have on $\cc^n \times \hh_n$ 
the theta function with \emph{characteristic} $\left[{a \atop b} \right]$ given by
\[ \theta \left[ {a \atop b } \right](z,\tau) = \sum_{m \in \zz^n} e^{\pi i
{}^t(m+a) \tau (m+a) + 2 \pi i {}^t (m+a)(z+b) } \, . \]
The choice $a=b=0$ gives the Riemann theta function, and it follows
from the definition that $\theta[{ a \atop b}](-z,\tau)=e^{4\pi i {}^t a b}
\theta[{a \atop b}](z,\tau)$. We call $[{a \atop b}]$ an even or odd theta
characteristic depending on whether 
$\theta \left[ {a \atop b } \right](z,\tau)$ is an even or odd function of $z$.
If $n=2$, there are ten even theta characteristics, and six odd ones. 
The product 
$\prod_{\varepsilon \, \textrm{even} } \theta[\varepsilon](0,\tau)^2$ 
of Thetanullwerte is a modular form of weight $10$ and level $1$ 
and can be related to the discriminant of
a hyperelliptic equation.
\begin{thm} For $(z,\tau)$ in $\cc^2 \times \hh_2$ with $\theta(z,\tau)=0$, 
the formula
\[ \eta(z,\tau)^3 = \pm \pi^{12} \prod_{\varepsilon \, \textrm{\emph{even}} }
\theta[\varepsilon](0,\tau)^2 \cdot \theta(3z,\tau) \]
holds.
\end{thm}
\begin{proof} It suffices to prove the formula for $\tau$ corresponding to an
indecomposable abelian surface. We write
$A=\cc^2/(\zz^2+\tau \zz^2)$ and $\Theta = \div \theta$ and assume that
$\Theta$ is 
irreducible. By Theorem \ref{main}, the section 
$\eta(z)$ is a theta function of order $3$ on $\Theta$, 
and by Example \ref{wronskiangenustwo} or Theorem \ref{jaco} it has
zeroes exactly at $\Theta \cap A[2]$, the Weierstrass points of $\Theta$, all of
multiplicity $1$. On the other hand, the function $\theta(3z)$ gives rise
to a global section of
$O_A(\Theta)^{\otimes 9}$, i.e. is a theta function of order $9$ on
$A$, and has zeroes on $\Theta$
exactly at $\Theta \cap A[2]$, with multiplicity $3$. It follows that $\eta(z)^3
= c \cdot \theta(3z)$ on the zero locus of $\theta(z)$ in $\cc^2$, 
where $c$ is a constant only depending on $\tau$. In order to compute $c$, we
recall from Remark \ref{formula} that $\omega_\zeta(z)^2 \eta(z) = F(z)^3$, with $\zeta$ the basis of
$H^0(\Theta,K_\Theta)$ given by $(dz_1,dz_2)$, so that
\begin{equation} \label{definingc}
c = \frac{ F(z)^9 }{ \omega_\zeta(z)^6  \theta(3z) } 
\end{equation}
for $z$ on $\Theta \setminus \left( \Theta \cap A[2] \right)$. We find $c$ by
letting $z$ approach a point $Q$ of $\Theta \cap A[2]$, along $\Theta$, 
and computing Taylor
expansions of the numerator and denominator in (\ref{definingc}). 
Note that the leading coefficient of a Taylor expansion of
$F(z)$ around $Q$ is the same as the leading coefficient of a Taylor expansion
of $\theta(2z)|_\Theta$ around $Q$. We start with the standard euclidean
coordinates $z_1,z_2$, but now translated suitably so as to have them both
vanish at $Q$ on $A$. According to \cite{grant}, formula
(1.6) there exist a constant $b_3$ and an invertible $2$-by-$2$ matrix $\mu$ such
that in the coordinates $(u_1,u_2) = (z_1,z_2) {}^t \mu$ 
one has a Taylor expansion
\begin{equation} \label{taylortheta}
   \theta(z) = \theta(\mu^{-1}u) = \gamma \, e^{G(u)} \left( u_1 +
\frac{1}{24}b_3u_1^3 - \frac{1}{12}u_2^3 + \textrm{higher order terms}
\right) \, , 
\end{equation}
with $\gamma$ some non-zero constant and with $G(u)$ some holomorphic function
that vanishes at $u=0$. This gives $u_2$ as a local
coordinate around $Q$ on $\Theta$, as well as an expansion
\begin{equation} \label{developu1}
u_1 = \frac{1}{12} u_2^3 + \textrm{higher order terms} 
\end{equation}
locally around $Q$ on $\Theta$. The way $u_1,u_2$ are obtained is as follows.
One may identify $(\Theta,Q)$ with a hyperelliptic curve $(X,\infty)$ given by a
hyperelliptic equation $y^2=f(x)$ with $f$ a monic separable 
polynomial of degree $5$. We
have then around $Q$ a local coordinate $t$ such that 
$x=t^{-2} + \textrm{h.o.t.}$ and $y=-t^{-5} + \textrm{h.o.t.}$ A computation
yields that 
\[ \int_\infty \frac{dx}{y} = \frac{2}{3} t^3 + \textrm{h.o.t.} \quad
, \quad \int_\infty \frac{xdx}{y} = 2t + \textrm{h.o.t.} \]
so that putting
\[ u_1 = \int_\infty \frac{dx}{y} \quad , \quad u_2 = \int_\infty
\frac{xdx}{y}  \]
gives the required relation (\ref{developu1}). 
The matrix $\mu$ corresponds then to a change of basis of holomorphic
differentials on $\Theta$ from
$\zeta=(dz_1,dz_2)$ to $\zeta'=(dx/y,xdx/y)=(du_1,du_2)$. From
(\ref{taylortheta}) and (\ref{developu1}) one computes
\[ \theta(2z)|_\Theta = \gamma \, e^{G(2u)} \left( 2u_1 - \frac{1}{12}(2u_2)^3 +
\textrm{h.o.t.} \right) = -\frac{1}{2} \gamma \, e^{G(2u)} \cdot u_2^3 + \textrm{h.o.t.} \]
and
\[ \theta(3z)|_\Theta = \gamma \, e^{G(3u)} \left( 3u_1 - \frac{1}{12}(3u_2)^3 +
\textrm{h.o.t.} \right) = -2 \gamma \, e^{G(3u)} \cdot u_2^3 + \textrm{h.o.t.}  \]
As to the Wronskian of $\zeta$ we have
\[ \omega_\zeta = (\det \mu)^{-1} \omega_{\zeta'} \, . \]
Writing out $\omega_{\zeta'}$ with respect to $u_2$ gives
\[ \omega_{\zeta'}(u_2) = \det \left( \begin{array}{cc} u'_1 & 1 
\\ u''_1 & 0 \end{array} \right) (du_2)^{\otimes 3} =
-u''_1 \, (du_2)^{\otimes 3} = 
\left( -\frac{1}{2} u_2 + \textrm{h.o.t.} \right) \, (du_2)^{\otimes 3} \, . \]
We find
\[ c = \frac{ F(z)^9 }{ \omega_\zeta(z)^6  \theta(3z) } =
\frac{ 2^{-9} \gamma^9 }{ 2^{-6} (\det \mu)^{-6} \cdot 2 \gamma}
= 2^{-4} (\det \mu)^6 \gamma^8 \, . \]
By \cite{grant}, Theorem 2.11 one has
\[ \gamma^8 = \pm 2^4 \pi^{12} (\det \mu)^{-6} \prod_{\varepsilon \, \mathrm{even}}
\theta[\varepsilon](0,\tau)^2 \, . \]
Substituting this in our formula for $c$ we get the result.
\end{proof}

\noindent \textbf{Acknowledgements} I am grateful to 
Prof. Ciro Ciliberto for his suggestions that led to 
Section~\ref{interpretation} of this paper. 
Furthermore I thank  
Samuel Grushevsky and Riccardo Salvati Manni for fruitful 
discussions and for several improvements upon earlier versions of Proposition
\ref{decomposables}.

\noindent Address of the author: \\  \\
Robin de Jong \\
Mathematical Institute \\
University of Leiden \\
PO Box 9512 \\
2300 RA Leiden \\
The Netherlands \\
Email: \verb+rdejong@math.leidenuniv.nl+

\end{document}